\input amstex
\magnification \magstep 1
\baselineskip=18truept

% Definitions

\def\pee#1{\hbox{$ {\bold P} ^ {#1}$}}

\def\pew#1#2{\hbox{$ {\bold P} ^ {#1}_{#2}$}}
\def \bign {\bigskip\noindent}
\def \medn {\medskip\noindent}
\def \parn {\par\noindent}
\def \tab#1{\kern #1 truein}
\def\qed{\hbox{$\clubsuit$}}

%fonts
\loadeufm  %usage is \frak
\loadmsbm  % usage is \Bbb

\def\OP#1{\hbox{${\Cal O}_{{\bold P}^{#1}}$}}

\def\E{\hbox{${\Cal E}$}}
\def\F{\hbox{${\Cal F}$}}
\def\G{\hbox{${\Cal G}$}}

\def\M{\hbox{${\Cal M}$}}

\def\Hom{\hbox{\rm Hom}}
\def\sHom{\hbox{\underbar{\rm{Hom}}}}

\def\I#1{\hbox{$\Cal I_{#1}$}}  %the ideal sheaf of a subscheme
\def\O#1{\hbox{$\Cal O_{#1}$}}  %the structure sheaf of a scheme 
  %the normal sheaf of a subscheme
% \def\T{{\Cal T}_{{\bf P}^n}}
% End definitions  **********************************************

\font\bigf=cmr10 scaled\magstep2
\vskip .5 truein
{\bigf \centerline
{Mathematical Instantons In Characteristic Two.}
}
\vskip .5 truein
{\bf \centerline
{ A.P.Rao}
}
\centerline{Department of Mathematics,}
\centerline{University of Missouri -- St\.Louis,}
\centerline{St\.Louis, MO 63121.}
\centerline{rao\@arch.umsl.edu}
%\vskip 1 truein
\font\smallf=cmsl10
{\smallf \underbar{Abstract}. On \pee3, we show that mathematical instantons 
in 
characteristic two are unobstructed.
We produce upper bounds for the dimension of the moduli space of stable rank
two bundles on \pee3 in characteristic two. In cases where there is a phenomenon
of good reduction modulo two, these give similar results in characteristic zero. We also give an example of a non-reduced component of the moduli space 
in characteristic two.
}
\vskip 1 truein

\parn
\underbar{{\bf Introduction}}
\medskip The study of mathematical instantons on projective three space has
been pursued partly because of the Atiyah-Ward-Drinfeld-Manin theorem which
showed that the solutions of the self-dual Yang-Mills equations on $S^4$
 could be described in
algebraic terms as particular cases of mathematical instanton bundles on
\pew 3{\Bbb C}. As 
a consequence, many workers have studied these bundles (and their generalizations
to \pee{2n+1})  over the complex numbers
and are still studying issues like the smoothness and irreducibility of the
moduli spaces of these bundles. In this paper, we would like to discuss
mathematical instanton bundles on \pee3 defined over
fields of any characteristic. 
The question of whether such bundles are unobstructed  is still
unknown in general, though it has been verified in some special cases. In [N-T],
unobstructedness is proved for all mathematical instantons with a section
in degree one. In [R-2], it is proved for those with a jumping line of maximal
order.  
\par
Over a field of characteristic
two, we will find (Theorem 2.4) that the unobstructedness of all mathematical instantons on \pee3
is  extremely easy to see. In fact, these are the only unobstructed stable 
rank two bundles on \pee3 (with $c_1=0$) in this characteristic. We also
 find 
a simple example (2.6) of a non-reduced component of the moduli space in characteristic two. We expect that most components will have this property of non-reducedness. 
 \par 
As a consequence, we show that in characteristic zero, a mathematical 
instanton is unobstructed if some pull-back of it by an automorphism of \pee3
reduces modulo 2 to a mathematical instanton (Theorem 3.2). However, this result does not answer the
problem for all mathematical instantons in characteristic zero. In fact,  there
are examples, in characteristic zero,  of such bundles for which no pull-back
by an automorphism of \pee3 reduces modulo two
to a mathematical instanton (Example 3.7).
\par
These computations in characteristic two also allow us to
bound the dimension of each component of the moduli space of stable rank two
bundles on \pee3. The bound has order $n^2$ where $n$ is the normalized second
Chern class (Corollary 2.8). Once again this bound has a consequence in characteristic zero.
Specifically, if $N$ is a component of the moduli space
which contains one bundle that reduces to a stable or semi-stable bundle in
characteristic two, then the dimension of $N$  is bounded above by a bound
of the order $n^2$ (Theorem 3.8).
\par
The first section contains a review of facts about bundles on \pee3\  in any characteristic, and  includes 
a definition of mathematical instantons (1.5). In section two, some calculations
in characteristic two are made. These arise from the relationship between
the second symmetric power in characteristic two and Frobenius pull-backs.
I would like to thank V. Mehta for pointing out to me  this relationship. In the last section,
applications to characteristic zero are made. 
\bign
{\bf \S 1}
\parn
We first review some elementary facts about vector bundles over a projective 
space defined over an arbitrary field $k$. Many of the results in the literature
are discussed for an algebraically closed field, and in the case of mathematical
instanton bundles, even over the complex numbers. We observe that most of these
conditions can be relaxed.

\medskip
	Consider $\pew nk$, projective $n$-space defined over a field $k$. Let
$\bar k$ be the algebraic closure  of $k$. The notion of a (geometric) vector
bundle over $\pew nk$ and the notion of a locally free sheaf on $\pew nk$
are equivalent ([H-1], II,5.18). Let \E\ be a vector bundle of rank $r$ defined on
$\pew nk$. Let $\bar {\E} = \E \otimes _k \bar k$ be its pull-back to \pew n{\bar k}. 
\medn
1.1) The square
$$\CD  \pew n{\bar k} @>>> \pew n k \\
         @VVV               @VVV   \\
        {\text Spec} \bar k @>>> {\text Spec}k \endCD $$
is Cartesian with $k @>>> \bar k$ faithfully flat. Hence 
$$ H^i( \pew n{\bar k}, \bar{\E}) = H^i(\pew nk, \E) \otimes_k \bar k$$
and
$$  H^i( \pew n{\bar k}, \bar{\E}) = 0  \iff  H^i(\pew nk, \E)  = 0.$$
\medn
1.2) The integers $c_i(\E)$ can be defined as $c_i(\bar{\E})$. However since
Pic(\pew nk)$\cong \Bbb Z$, there is an isomorphism $\wedge^r \E \cong \O {\pew nk}
(c_1)$ which is defined over $k$.
\medn
1.3) Horrock's Theorem states that \E\ is isomorphic over $k$ to a sum of
line bundles if and only if $ H^i_*(\pew nk, \E) = 0$ for all $i$ between $1$
and $n-1$.
This is valid over any field. For example, consider the proof given in [O-S-S], 
which uses the complex numbers as the base field.
Upon reading the proof, we see only one place where the 
argument does not work for arbitrary $k$. This is in the proof of Grothendieck's 
theorem, an auxiliary result needed in the proof. In this part, a section $s \in  H^0(\pew 1k, \E(k_0))$ is chosen, where $k_0$ is the
least integer for which a nonzero section $s$ can be chosen and the claim  is
made that this
section is nowhere vanishing. For us, in our context where $k$ is arbitrary, this should mean that $s$ has no zeros over 
$\bar k$. Indeed this is true, for if $s$ has a zero in \pew 1 {\bar k}, we would conclude that
$s$ comes from  a section of $ H^0(\pew 1{\bar k}, \E(k_0-1))$. But then
$H^0(\pew 1k, \E(k_0-1))$ itself is non-zero by (1.1) contradicting our choices.
Thus we still get $$ 0 @>>> \O {\pew 1k} @>s>> \E(k_0) @>>> \F @>>> 0 $$
where \F\ is a bundle on \pew 1k, and the proof in [O-S-S] continues without 
change.
\medn
1.4) Therefore the results of [B-H] and [R-1] on the construction of monads for a bundle \E\  
of rank two on \pew 3k are valid in any characteristic. Let \E\ have first Chern 
class $c_1$. Then there is an isomorphism $\E^\vee \cong \E(-c_1)$ which is defined
over $k$. $M = H^1_*(\pew 3k, \E)$ is a finite length module over $S=k[X_0,X_1,X_2,
X_3]$ . Let $L_0 @>>> M$ be a surjective homomorphism where $L_0$ is a sum of
graded twists of $S$, picking out a set of minimal generators of $M$. Then
there is a monad
$$ 0 @>>> \tilde L_0^\vee(c_1) @>\beta >> \tilde L_1 @>\alpha >> \tilde L_0 @>>> 0$$
where $L_1$ is also a sum of twists of $S$ and $\alpha, \beta$ are matrices of
homogeneous polynomials in $S$. Furthermore, there is an isomorphism
$H: L_1 \cong L_1^\vee(c_1)$ with $H$ a matrix of homogeneous polynomials in $S$
such that $H\beta = 
\alpha^\vee.$ This gives an isomorphism between the monad and the dual monad which lifts
the isomorphism between $\E^\vee$ and $\E(-c_1)$. If $\Hom(L_0, L_1) = 0$, then
this $H$ is unique and can be chosen so that $H^\vee = -H$. 
\par
The result in [R-1]
says that we may take $\alpha$ as a minimal presentation of the $S$-module
$M$. Furthermore if $L_2 @>>> L_1 @>>> L_0 @>>> M @>>> 0$ is part of a minimal resolution
for $M$, then the map $ L_0^\vee(c_1) @>\beta >> L_1$ will be a direct summand of the map $L_2 @>>> L_1 $ [H-R, 3.2].  
\medn
1.5) We give the definition of a mathematical instanton of rank two on 
\pew 3 k . 
The general definition due to Okonek and Spindler [O-S] of a mathematical 
instanton on \pee{2n+1} has conditions on natural cohomology and trivial splittings on the general line. We
will relax this part of their definition for our case of rank two on \pee3.
\medn
\underbar{Definition}: An indecomposable  rank two bundle \E\ on \pew 3k with $c_1 = 0, c_2 =n$ is called a mathematical
instanton if $H^1(\pew 3k, \E(-2)) = 0$. 
\par
Recall that for a rank two bundle with $c_1 = 0$ in characteristic
zero, the condition that it splits trivially on the general line is equivalent to the condition of semi-stability. Likewise, natural cohomology in the range 
$-3$ to 0 also implies stability (if $H^1(\E) \neq 0$.)  We will not assume  
{\it a priori} such conditions. However, stability will follow below.  
In arbitrary characteristic, we do not know if the splitting type of  
a bundle as defined above can be non-trivial on the general line. 
The following theorem is well known. We include
a proof because, for example, the proof in [B-H] uses trivial splitting on
the general line and the corresponding proof in [O-S] uses their condition
of natural cohomology. 
\medn
\underbar {Theorem}: Let \E\ be a mathematical instanton on \pew 3k with $c_1 = 0, c_2 =n$. Then
\par
a) $H^1(\pew 3k, \E(-k)) = 0$ for $k \geq 2$.
\par
b) $M=  H^1_*(\pew 3k, \E)$ has all its minimal generators in degree -1.
\par
c) $H^0(\pew 3k, \E) = 0$ (hence \E\ is stable.)
\par
d) $n = h^1(\pew 3k, \E(-1)) > 0$.
\par
e) \E\ is the homology of a minimal monad
 $$ 0 @>>> n\O {\pew 3k}(-1) @>\beta>> (2n+2)\O {\pew 3k} @>\alpha >> n\O {\pew 3k}(+1)
@>>> 0.$$
\par
f) Conversely, any bundle which is the homology of such a monad as in (e) is 
a mathematical instanton bundle with $c_1 = 0,c_2 =n$.
\medn
\underbar{Proof}: 
Since most of the statements are about dimensions of
cohomology, by (1.1) we will assume that $k$ is algebraically closed so that 
geometric constructions work as usual. To prove (a), we will show
that if $c_1(\E) = 0$ or $-1$ and $H^1(\pew 3k, \E(-m)) = 0$ for some $m>0$, then 
$H^1(\pew 3k, \E(-k)) = 0$ for all $k \geq m$. For suppose 
$H^1(\pew 3k, \E(-k)) \neq 0$ with $k>  m$ and let $k$ be the least such. Let
$H$ be a general hyperplane in \pew 3k, and consider
$$ 0 @>>> \E(-k) @>H>> \E(-k+1) @>>> \E_H (-k+1) @>>> 0. \tag * $$
Then $\E_H (-k+1)$ clearly gets at least one global section $t$. Consider
$$  0 @>>> \E(-k+1) @>H>> \E(-k+2) @>>> \E_H (-k+2) @>>> 0.$$
Since $H^1(\E(-k+1)) = 0$, the multiples of $t$ in $H^0(\E_H (-k+2))$ arise from
 sections of $\E(-k+2) $. So $\E(-k+2) $ has at least three global sections. 
The sections of \E\ in degrees less than or 
equal to 0 (if they exist) are all multiples of a single section $s$ and this
section $s$ is the unique section in degrees $\leq 0$ whose zero-scheme has codimension
two. Since $-k+2 \leq 0$, there
is a single section, say $s$,  of
$\E(-l)$ for some $l >0$. This induces a non-zero section $s'$ of $\E_H (-l)$.
If $H$ is chosen generally,
the zero-scheme of $s'$ in $H$ has codimension two. Hence all sections of $\E_H$ in degrees
$\leq 0$ must be multiples of $s'$. It is evident from the long exact sequence of
cohomology of (*) that the section $t$ of
$\E_H(-k+1)$ is not a multiple of $s'$. This is a contradiction. Hence the result is 
proved.
\par
To prove (b), we study the minimal monad of \E. The monad will be $0 @>>> \tilde
L_0^\vee @>>> \tilde L_1 @>>> \tilde L_0 @>>> 0$ as before. If $M$ has any 
generators in degrees 0 or 1 etc, it means that $L_0$ has summands like $S(0)$
or $S(-1)$ etc. By the minimality of the monad, $L_1$ must have summands like
$S(-1)$ or $S(-2)$ etc. But $L_1^\vee \cong L_1$, hence $L_1$ has summands like
$S(1)$ or $S(2)$ etc. Such summands must map to zero in $L_0$ by degree considerations
and minimality. However, the map $L_1 @>>> L_0$ is a minimal presentation of $M$,
hence no summands can map to zero. Hence all generators of $M$ are in degree $-1$.
\par
To prove (c), we now know that $L_0$ is a sum of $S(1)$'s. Hence $L_1$ cannot have 
any summands like $S(1)$ or $S(2)$ etc. Since $L_1$ is self dual, $L_1$ can contain
only $S(0)$'s. Hence the minimal resolution of $M$ which looks like
$ @>>> L_2 @>>> L_1 @>>> L_0 @>>> M @>>> 0$ must have $L_2$ without any $S(0)$
or $S(1)$ etc.  Now $L_2 = L_0^\vee \oplus L_2'$ and there is a surjection
of $\tilde  L_2' @>>> \E$ which induces a surjection of global sections in all
twists. This shows that \E\ has no sections in degree $\leq 0$.
\par
(d) is merely a Riemann-Roch computation now. $n>0$ because  for example, 
if $n=0$, then we get $M=0$ by (b). Hence by Serre duality, $H^2_*(\E) = 0$ as well. 
So \E\ is decomposable by Horrock's theorem. Therefore $n>0$.
\par Of course by now, (e) has been demonstrated. 
(f) is quite obvious from the display of the monad.\qed

\medn
1.6)  
The coarse moduli scheme $\M_{\pew 3k}(c_1,c_2)$ of stable rank 2 bundles on \pew 3 k with Chern classes 
$c_1, c_2$ exists and is a 
$k$-scheme ([M], Theorem 5.6). 
The fact that this moduli scheme is  quasi-projective has been proved
by Maruyama and also discussed in [H-2]. This scheme behaves well under base
change ([M], Remark 5.9) so that for example if $\bar k$ is the algebraic closure of $k$, then
$\M_{\pew 3 {\bar k}}(c_1,c_2) \cong \M_{\pew 3k}(c_1,c_2) \times_{{\text{Spec}k}}
{\text {Spec} \bar k}$.
If \E\ is a bundle on \pew 3k giving a $k$-valued point on $\M_{\pew 3k}(c_1,c_2)$,
then  $ \bar {\E} =\E \otimes_k \bar k$ on \pew 3 {\bar k} gives
a geometric point. The Zariski tangent space at this geometric point is given by
$H^1(\pew 3{\bar k}, \sHom(\bar{ \E}, \bar{\E}))$ and if $H^2(\pew 3{\bar k}, \sHom( 
\bar{\E}, \bar{\E}))=0$, the moduli
space is smooth at this point, with  dimension equal to $h^1(\pew 3{\bar k}, 
\sHom(\bar {\E}, \bar{\E}))$. In this case, we shall say that \E\ is unobstructed. 
Of course, the dimensions of these vector spaces can be computed using \E\ over $k$.

\medn
1.7)
The spectrum of a semi-stable rank two bundle on \pew 3k has been defined in characteristic
zero in [B-E] and in arbitrary characteristic in [H-3]. Since the spectrum determines and is determined by dimensions of $h^1(\pew 3k, \E(l))$, we do not need to
assume that $k$ is algebraically closed. Let us recall the spectrum as found
in [H-3].
Let \E\ be a semi-stable rank two bundle on \pew 3k with $c_1 = 0$ or $-1$ and $c_2 = n$.
There is a unique set of $n$ integers $k_1, k_2, \dots k_n$ called the spectrum
of \E\ with the following properties: let $\Cal H$ denote the sheaf $\oplus \O
{\pee1}(k_i)$ on \pee1.
\parn
a) $h^1(\pew 3k, \E(l)) = h^0(\pee1, \Cal H(l+1))$ for $l \leq -1$.
\parn
b) $\{-k_i\} = \{k_i\}$ if $c_1 = 0$ and
\parn
\ \ $\{-k_i\} = \{k_i +1\}$ if $c_1 = -1$.
\parn
c) The spectrum is connected, except possibly for a gap at $0$. If \E\ is stable,
then the spectrum is connected.
\parn
d) An integer $l$ may appear more than once in the spectrum. If $l \leq -2$ and
$l$ appears exactly once in the spectrum, then any smaller integer can occur at
most once in the spectrum. If \E\ is stable and $c_1 = 0$, we can say the same for
$l \leq -1$. 
\par
(d) was  proved in [H-3 Proposition 5.1] using a characteristic zero hypotheses. 
However, this hypotheses was really needed  only to prove a stronger statement about unstable
planes and as was pointed out in [H-R], for the proof of (d), characteristic zero
is not required.
\medn
1.8) When $c_1=0$, the condition that $H^1(\pew 3k, \E(-2)) = 0$ is equivalent 
to the condition that the spectrum of \E\ (with $c_2 = n$) consists of $n$ zeroes.
Likewise, let \E\ be a rank two bundle with $c_1 = -1, c_2 = n$ and with 
$H^1(\pew 3k, \E(-2)) = 0$. In our proof of Theorem 1.5, we showed that
$H^1(\pew 3k, \E(-k)) =0$ for all $k \geq 2$. It is easy to prove by the same techniques that, 
in the minimal monad for \E, the term $L_1$ has only $S(0)$'s and $S(-1)$'s,
hence $L_2$ can contain only terms $S(a)$ with $a \leq -1$. Therefore, \E\ is  
stable.
\E\ has spectrum consisting of $\frac n2$ 0's and   $\frac n2$ $-1$'s.

\bign
{\bf \S 2}
\parn
Let $\pi:X @>>> Y$ be a morphism of schemes defined over a field $k$ of characteristic
$p$ different from zero. We can define Frobenius automorphisms $F$ of $X$ and $Y$
induced by the Frobenius homomorphism $a \longmapsto a^p$ on affine rings.
Then the square
$$ \CD X @>F>> X\\
      @V\pi VV  @V\pi VV\\
       Y @>F>> Y \endCD $$
commutes.
\medn
2.1)\underbar{Lemma}: Let $F: \pew nk @>>> \pew nk$ be the Frobenius morphism in
characteristic $p$. Then
\parn
(i) $F_*[\O {\pew nk}(l)]  \cong \oplus_{i\geq -\frac lp}  a_i \O {\pew nk} (-i)$ where $a_i$
is the number of monomials $X_0^{b_0}X_1^{b_1}\dots X_n^{b_n}$
of degree $l+pi$ with each exponent $b_j < p$.
\parn
(ii)  When  $k$ has characteristic two, 
$$ F_*[\O {\pew nk}(l)] \cong \bigoplus_{\frac {n+1-l}2 \geq i\geq -\frac l2 } \binom {n+1}{l+2i}
\O {\pew nk} (-i)$$.
\parn
\underbar{Proof}: For $1 \leq i \leq n-1$, $$ \align  H^i(\pew nk, F_*[\O {\pew n
k}(l)] \otimes \O {\pew nk}(m)) &= H^i(\pew nk, \O {\pew nk}(l) \otimes F^*[\O {\pew nk}
(m)])\\ & = H^i(\pew nk,  \O {\pew nk}(l+pm)) = 0 \endalign $$ for all $m$.
Hence by Horrock's
theorem, the
bundle $ F_*[\O {\pew nk}(l)]$ is a sum of line bundles.
The number of summands in this bundle of the form $\O {\pew nk} (-i)$ can be 
computed
by finding the dimension of $H^1(\pew nk,  F_*[\O {\pew nk}(l)]\otimes \Omega^1_{\pew
nk}
(i))$ which is just $h^1(\pew nk, F^*[\Omega^1_{\pew nk}](l+pi))$. We have the 
sequence
$$ 0 @>>>  F^*[\Omega^1_{\pew nk}](l+pi) @>>> (n+1)\O {\pew nk}(l+p(i-1)) @>{\bmatrix
X_0^p,\dots,X_n^p \endbmatrix }>> \O {\pew nk}(l+pi) @>>> 0.$$ The lemma follows
now from this sequence. \qed

\medn
2.2)\underbar{Proposition}:
Let \E\ be a rank two bundle on \pew nk with first Chern class $c_1$, where
$k$ has characteristic $p$. Let $F$ be the Frobenius mapping on  \pew nk. Then we have exact sequences
$$ 0 @>>> F^* \E @>>> S_p(\E) @>>> S_{p-2}(\E)\otimes \O {\pew nk}(c_1) @>>> 0$$
and 
$$   0 @>>> \O {\pew nk} @>>> \sHom(\E,\E) @>>> S_2\E \otimes \O {\pew nk}(-c_1)
@>>> 0. $$
\parn
\underbar{Proof}:
Consider the commuting square (not Cartesian)
$$ \CD  \bold P(\E) @>F>> \bold P(\E) \\
         @V\pi VV        @V\pi VV \\
         \pew nk @>F>>    \pew nk .\endCD $$
Let $\O {\pi}(1)$ denote the tautological line quotient bundle on $\bold P(\E)$.
So we have
$$ 0 @>>> \wedge^2(\pi^* \E) \otimes \O {\pi}(-1) @>>> \pi^*\E @>>> \O {\pi}(1)
@>>> 0.$$
Since $\wedge^2(\pi^* \E) \cong \pi^*\O {\pew nk}(c_1)$, after applying $F^*$ we get
$$ 0 @>>> F^*\pi^*\O {\pew nk}(c_1) \otimes F^*\O {\pi}(-1) @>>> F^*\pi^*\E
@>>> F^*\O {\pi}(1) @>>> 0, $$
hence
$$  0 @>>> \pi^*\O {\pew nk}(pc_1)\otimes \O {\pi}(-p) @>>> \pi^*F^*\E
@>>> \O {\pi}(p) @>>> 0.$$
Applying $\pi_*$, we get (since $\pi_*(\O {\pi}(-p))=0$)
$$ 0 @>>> F^*\E @>>> \pi_* \O {\pi}(p) @>>> \O {\pew nk}(pc_1)\otimes R^1\pi_*\O {\pi}(-p) @>>> 0.$$
Now $R^1\pi_*\O{\pi}(-p) \cong [\pi_* \O{\pi}(p-2)]^\vee \otimes (\wedge^2\E)^\vee$ ([H-1], III,8.4),
hence we get
$$ 0 @>>> F^* \E @>>> S_p(\E) @>>> [S_{p-2}\E]^\vee \otimes \O {\pew nk}(pc_1-c_1) @>>> 0. $$

Now in characteristic $p$, $ [S_{p-2}\E]^\vee \cong S_{p-2}(\E^\vee) \cong
S_{p-2}(\E)\otimes \O {\pew nk}(-(p-2)c_1)$
hence we end with
$$ 0 @>>> F^* \E @>>> S_p(\E) @>>> S_{p-2}(\E)\otimes \O {\pew nk}(c_1) @>>> 0.$$\medskip
For the other part, observe that $\sHom(\E,\E) \cong \E^*\otimes \E \cong \E\otimes
\E\otimes \O {\pew nk}(-c_1)$, hence we have the sequence
$$   0 @>>> \O {\pew nk} @>>> \sHom(\E,\E) @>>> S_2\E \otimes \O {\pew nk}(-c_1)
@>>> 0$$ (obtained for example by the push down of the tautological sequence
on $\bold P(\E)$ tensored by $\O {\pi}(1)$.)
\qed
\medn
2.3)\underbar{Corollary}: Let \E\ be a rank two bundle on \pew 3k, where $k$ has 
characteristic two. 
\parn
i) Let  $c_1 = 0$ and let $m\geq -4$.
\par
 If $m$ is even, then $h^2(\pew 3k, \sHom(\E,\E)(m)) =  h^1(\pew 3k, \E(-2-
\frac m2)) + 6 h^1(\pew 3k, \E(-3-\frac m2)) +  h^1(\pew 3k, \E(-4-\frac m2)). $ 
\par
 If $m$ is odd, then $h^2(\pew 3k, \sHom(\E,\E)(m)) =  4 h^1(\pew 3k, \E(-\frac {3+m}
2 -1) + 4  h^1(\pew 3k, \E(-\frac {3+m}2 -2))$.
\parn
ii) Let  $c_1 = -1$ and let $m\geq -4$.
\par
 If $m$ is even, then $h^2(\pew 3k, \sHom(\E,\E)(m)) =  4 h^1(\pew 3k, \E(-2 -
\frac m2)) +4 h^1(\pew 3k, \E(-3 -\frac m2)).$
\par
 If $m$ is odd, then $h^2(\pew 3k, \sHom(\E,\E)(m)) = h^1(\pew 3k, \E(-\frac {m+3}2)
+ 6  h^1(\pew 3k, \E(-\frac {m+3}2 -1)) +  h^1(\pew 3k, \E(-\frac {m+3}2 - 2)).$

\underbar{Proof}: When $m\geq -4$, $$ \align h^2(\pew 3k, \sHom(\E,\E)(m)) &= 
h^2(\pew 3k, S_2\E \otimes \O {\pew 3k}(m-c_1))\\ &=
h^2(\pew 3k, F^*[\E]\otimes \O {\pew 3k}(m-c_1)) \\
&= h^1(\pew 3k, F^*(\E^\vee) \otimes \O {\pew 3k}(c_1-m-4)) \\
&= h^1(\pew 3k, F^*(\E)  \otimes \O {\pew 3k}(-2c_1 +c_1 -m -4)) \\
&= h^1(\pew 3k, \E \otimes F_*[\O {\pew 3k}(-c_1-m-4)]). \endalign $$
Now if $-c_1-m-4 = 2t$, then $F_*[\O {\pew 3k}(-c_1-m-4)] = F_*F^*[\O {\pew 3k}
(t)] = \O {\pew 3k}(t) \otimes F_*[\O {\pew 3k}] = \O {\pew 3k}(t) \otimes [
\O {\pew 3k}\oplus 6\O {\pew 3k}(-1) \oplus \ {\pew 3k}(-2)]$.
If  $-c_1-m-4 = 2t-1$, then 
$F_*[\O {\pew 3k}(-c_1-m-4)] =  F_*[F^*[\O {\pew 3k}(t)] \otimes \O {\pew 3k}(-1)]
= \O {\pew 3k}(t) \otimes F_*[ \O {\pew 3k}(-1)] 
=  \O {\pew 3k}(t) \otimes [ 4  \O {\pew 3k}(-1) \oplus 4 \O {\pew 3k}(-2)].$
\qed

\medn
2.4)\underbar{Theorem}: If $k$ is a field of characteristic two, then a stable rank
two bundle \E\ with $c_1=0$ is unobstructed if and only if \E\ is a mathematical
instanton on \pew 3k. 
\parn
\underbar{Proof}: 
$h^2(\pew 3k, \sHom(\E,\E)) = h^1(\E(-2)) + 6 h^1(\E(-3)) + h^1(\E(-4))$ by the 
Corollary above. We saw earlier (Theorem 1.5 (a))
 that $h^1(\E(-2)) = 0$ implies that the other
terms are also zero. \qed
\medn
2.5)\underbar{Theorem}: If $k$ is a field of characteristic two, then a stable rank
two bundle \E\ with $c_1=-1$ is unobstructed if and only if $h^1(\pew 3k, \E(-2))
=0$.
\parn
\underbar{Proof}: In this case, $h^2(\pew 3k, \sHom(\E,\E))=
4 h^1(\E(-2)) + 4  h^1(\E(-3))$. Again, if $h^1(\E(-2)) =0$ then so is
$h^1(\E(-3))$ (see 1.8). \qed

\medskip
With this situation, one expects that the moduli schemes in characteristic two
will be highly singular. Indeed, in the first example of a bundle family not
of mathematical instanton type, we find that the moduli scheme has a non-reduced
component. In contrast, at this time, very few examples of singular components
 of the moduli scheme are known in characteristic zero ([Ma], [A-O]).

\medn
2.6)\underbar{Example}: A non-reduced component of $\M(0,3)$ in characteristic two.
 \par Consider stable bundles on \pew 3k  with $c_1=0, c_2=3$ and spectrum
$-1,0,1$. The monad of such a bundle \E\ (regardless of characteristic)
has the form
 $$ 0 @>>> \O {\pew 3k}(-2) @>\alpha>> \O {\pew 3k}(-1) \oplus 2  \O {\pew 3k} \oplus
\O {\pew 3k}(1) @>\beta>> \O {\pew 3k}(2) @>>> 0. $$
Calling it $ 0 @>>> A @>\alpha >> B @>\beta >> C @>>> 0$, let
 \G\ be the kernel of the right hand map. \G\ has a ten dimensional family of
sections in degree 2, hence the set of all such monads is parametrized by a
quasi-projective variety $V$ of dimension $10 + \text{dim} Hom(B,C) = 54$. This
space maps onto the moduli space of all stable bundles of the type being 
considered. The group Aut($A) \times $Aut($B)\times $Aut($C$) acts on $V$ and
the orbit of a monad consists of monads for isomorphic stable bundles. 
The subgroup $k^*$, embedded diagonally, stabilizes a monad. On the other hand,
since $Hom(C,B) = Hom(B,A)=0$, results of Barth and Hulek [B-H] tell us that any
automorphism of a bundle \E\ is uniquely lifted to an automorphism of monads.
Since \E\ is stable, Aut(\E) consists of elements of $k^*$, hence the stabilizer
of a monad is exactly $k^*$. Thus we get a dimension of 54 less $(1+32 +1 -1)$
or a dimension of 21 for this component of the moduli space. 
\par
The mathematical instantons form smooth components of the moduli space, of dimension
21,  proved above in characteristic two. (See also, for example,  [L], for a
proof valid in any characteristic.) There are just two possible spectra for 
these Chern classes.
Hence, by reasons of dimension,  the bundles with spectrum
$-1,0,1$ give a distinct irreducible component of the moduli scheme.
\par
Now  these bundles \E\ (with spectrum
$-1,0,1$)  have  $h^1(\pew 3k, \E(-2)) = 1$,
Hence,  in characteristic two, by Corollary 2.3, 
$h^2(\pew 3k, \sHom(\E,\E)) = 1$. Since $h^1(\pew 3k, \sHom(\E,\E)) -
h^2(\pew 3k, \sHom(\E,\E)) = 8c_2 -3 = 21$, it follows that the Zariski tangent
space to each such  point \E\ on this component of the moduli scheme is 22-dimensional
and hence larger than the dimension of the component. Thus we get a non-reduced component. \qed

\medn
\medn
2.7)\underbar{Theorem}: 
Let \E\ be a semistable bundle of rank two with $c_1 = 0$ or $-1$ and $ c_2 = n$  on \pew 3k, in characteristic
2.
\parn
i) If $c_1 =0$, then  $h^2(\pew 3k, \sHom(\E,\E)) \leq (n-1)^2$.
\parn
ii) If $c_1 =0$ and \E\ is stable, then $h^2(\pew 3k, \sHom(\E,\E)) \leq (n-
2)^2 $.
\parn
iii) If $c_1 = -1$ (so \E\ is stable), $h^2(\pew 3k, \sHom(\E,\E)) \leq (n-2)^2$. 

\medn
\underbar{Proof}:
By the earlier discussion, we just need to bound $h^1(\pew 3k, (F^* \E)(-4))$
above.

First, let $c_1 = 0$. Using Lemma 2.1 , we see that $$ \align h^1(\pew 3k, (F^* \E)(-4))&= h^1(\pew 3k, F^*(\E(-2))) \\ &= h^1(\pew 3k, \E(-2) \otimes F_*(\O {\pew 3k} )
\\ &= h^1(\pew 3k, \E(-2)) + 6 h^1(\pew 3k, \E(-3)) + h^1(\pew 3k, \E(-4)) . \endalign  $$ 
Now \E\ has a spectrum of $n$ integers and let the
positive integers in the spectrum consist of $a_1$ ones, $a_2$ twos, \dots,
$a_r$ $r$'s with no $a_i$ equal to zero by the connectedness of the spectrum.
Let $\sum a_i =b$. Let $\Cal K = \oplus_{i>0} a_i \O {\pee1}(i)$. Then
$h^1(\pew 3k, \E(-m)) = h^0(\pee1, \Cal K(-m+1))$ for $m\geq 1$.
\par
Hence $$ \align  h^1(\pew 3k, \E(-2)) &= a_1 + 2a_2 + 3a_3 \dots + ra_r \\
         6h^1(\pew 3k, \E(-3)) &= 6(a_2 + 2a_3 \dots +(r-1)a_r \\
         h^1(\pew 3k, \E(-4)) &= a_3 + 2 a_4 \dots + (r-2) a_r, \endalign $$
and
$$ h^1(\pew 3k, (F^* \E)(-4)) = a_1 + 8(a_2 + 2a_3 + \dots + (r-1)a_r).$$
If $b$ is fixed, this sum will be maximized when each $a_i = 1$, giving
$$ h^1(\pew 3k, (F^* \E)(-4)) \leq 1 + 4(b-1)b = (2b-1)^2.$$ 
\parn
Now when \E\ is semistable, we know by symmetry of the spectrum that
$b \leq \frac n2$.  If in addition  \E\ is stable, then the spectrum is connected, hence  $b  \leq \frac {n-1}
2$. 
This gives (i) and (ii).
\parn
If $c_1 = -1$, we need to bound $4[h^1(\pew 3k, \E(-2)) + h^1(\pew 3k, \E(-3))]$which is equal to $4[a_1+3a_2+5a_3+\dots + (2r-1)a_r]$. For fixed $b =
\sum a_i$, this sum is maximized when all
$a_i$'s are 1, hence by $4b^2$. Now $b \leq \frac n2 -1$ since the spectrum is
connected and symmetric about $-1/2$. Hence the bound of (iii). \qed

\medn
2.8)\underbar{Corollary}:
In characteristic two, the moduli spaces of stable rank two bundles on \pew 3k
have each component bounded above in dimension: if $n$ is the normalized second
Chern class, the dimension is less than or equal to $n^2 + 4n +1$ for $c_1 = 0$ and less than or equal to $n^2 + 4n -1$ for  $c_1 = -1$.
\parn
\underbar{Proof}: If $c_1 = 0$ (respectively $-1$), then 
$h^1( \sHom(\E,\E)) - h^2(\sHom(\E,\E))$
equals $8n-3$ (respectively $8n-5$), for such a stable bundle \E. Now use the last theorem to bound $h^1
(\sHom(\E,\E))$. \qed

\medn
2.9) \underbar{Remarks}: It is quite likely that this upper bound is too 
coarse. In fact, one expects that non-reducedness contributes a part to
this bound for $h^1(\sHom(\E,\E))$.
  Known examples of large families of stable bundles have dimension 
much below this bound. Examples of Ellingsrud and Str\o mme give components
of $\M(0,2k-1)$ of dimension $3k^2+4k+1$, while Ein gives examples in
$\M(-1, 2k)$ of dimension $3k^2+7k+2$ (See [E] for both. $k \geq 2$ in these
examples.) These dimensions
are of the order $\frac 34 n^2$.
% Note that Ein includes the case $k=1$ in his examples, but his dimension is
%11 and not 12 in this case.
\bign
{\bf \S 3}
\parn
We will try to draw some conclusions about bundles in characteristic zero from
these results in characteristic two. So let \E\ be a rank two bundle on
\pew 3k,  defined
over a field $k$ of characteristic zero. 
% By the `Lefschetz principle', we may take $k$ to be a finitely 
% generated extension field of $\Bbb Q$. 
By the discussion
in the first section, \E\ is the homology of a minimal monad
$$ 0 @>>> \tilde L_0^\vee(c_1) @>\beta >> \tilde L_1 @>\alpha >> \tilde L_0 @>>> 0$$
where $\alpha, \beta$ are matrices of homogeneous polynomials in $S = k[X_0,X_1,
X_2,X_3]$. Let $A \subset k$ be a sub-integral domain of $k$ whose field of
fractions is $k$. Since the monad gives an equivalent monad if $\alpha, \beta$
are multiplied by nonzero elements of $k$, we may assume that 
 $\alpha, \beta$ are matrices of homogeneous polynomials in $A[X_0,X_1,X_2,X_3]$.
We will call this `a lift of the monad (and of the bundle \E) to $A$'. This lift
is of course by no means unique. Now let $\frak p$ be a prime ideal of $A$ such
that the residue field $k(\frak p)$ has characteristic two. For a given $A$,
such an ideal may not exist in $A$,
but there will always be an $A$ in $k$ for which such an ideal exists. 
Then, taking the lifted monad for \E, we may reduce it  modulo this ideal,
that is to say, we may apply $\otimes _A k(\frak p)$.
In general, we do not expect this to be a monad over the field $k(\frak p)$; for  example, the reduced 
$\beta$ may not be an inclusion of bundles. 
\medn
3.1)\underbar{Definition}: We will say that \E\ has  a good monad
reduction to characteristic two if there is a minimal monad for \E, as described above for some choice of $A$ and $\frak p$, such that the 
reduction modulo $\frak p$ is still a monad over the residue field.
\medn
\underbar{Remarks}:
Note that this definition, we are ending up with a bundle in characteristic two
which has the same monad type as the bundle we started with. So, for example,
if a mathematical instanton bundle has a good monad reduction to characteristic two, it's lift
to $A$  specializes to a mathematical instanton bundle in characteristic two.
\par
We can also assume that $(A, \frak p)$ is a discrete valuation ring. Indeed, by the 
`Lefschetz Principle', we may assume that $k$ (the field of definition of \E)
is finitely generated over $\Bbb Q$. If an $A$ has been found in $k$ with a
$\frak p$ giving good monad reduction, we can replace $A$ with its
integral closure $A'$ in $k$ and replace $\frak p$ with an over-ideal $\frak p'$ in 
$A'$, since reduction to $k(\frak p')$ is obtained by base change after first
reducing to $k(\frak p)$. Next, replacing $A'$ with its localization at 
$\frak p'$, we may assume that $(A, \frak p)$ is a local normal domain with
fraction field $k$. The condition of good reduction modulo $\frak p$ and
Nakayama's Lemma tell us that the lift of the monad to $A$ defines a monad
for a vector bundle on \pew 3 A, hence also at each localization $A_{\frak q}$
of $A$. Since $2 \in \frak p$, we can choose a height one 
prime sub-ideal $\frak q$  containing $2$,  
getting a discrete valuation ring $A_{\frak q}$ with the required properties.  
\medn
3.2)\underbar{Theorem}: Let \E\ be an mathematical instanton  bundle on \pew 3k, where
$k$ has characteristic zero.
\parn
1) If \E\ has a good monad reduction to characteristic two,
then \E\ is unobstructed. 
\parn
2) Suppose there is an element $\varphi \in GL(4,k)$ such that $\varphi^*\E$
has a good monad reduction to characteristic two. Then \E\ is unobstructed.
\parn
\underbar{Proof}: This is a standard upper semi-continuity argument. Let $A$, $\frak
p$ be as in the definition above. Then since modulo $\frak p$, the lift of the
 monad gives a monad over the residue field, there is an open set $U$ in Spec $A$
which contains $\frak p$, and such that on \pew 3 U, the lift of the monad is
a monad, ie. \E\ lifts to a bundle $\E_U$ on \pew 3 U. Since at the prime $\frak p$,
the restriction of $\E_U$ is a mathematical instanton, hence unobstructed (in
characteristic two), there is a perhaps smaller open set $V$ in $U$, where we may
take $V= \text{Spec} A_f$ for some $f\in A$, over which $H^2(\pew 3 V, \sHom(
\E_V,\E_V)) = 0$, hence also $H^2(\pew 3 k,  \sHom(\E, \E)) = 0$. 
\par
The second result follows since \E\ and $\varphi^*\E$ have appropriate cohomology groups
of the same dimension. \qed

\medn
\underbar{Examples}: 
\parn
3.3) On \pew 3{\Bbb Q}, consider the mathematical instanton bundle \E\ with $c_1 =
0, c_2 = 1$ and minimal monad
 $$ 0 @>>> \OP3 (-1) @>\beta >> 4 \OP3 @>\alpha >> \OP3 (1) @>>> 0, $$
where $\beta = \bmatrix -2X_1 \\ X_0\\-2X_3\\X_2 \endbmatrix, \alpha =
\bmatrix X_0 & 2X_1 & X_2 & 2X_3 \endbmatrix.$ Certainly, this monad (defined
over $\Bbb Z$) has bad reduction
modulo 2. However, \E\ has an equivalent monad which has good reduction, for we
see that there is an isomorphism of monads over $\Bbb Q$ given by
$$ \CD
0  @>>> \OP3 (-1) @>\beta >> 4 \OP3 @>\alpha >> \OP3 (1) @>>> 0 \\
@.        @V 2 VV             @V \gamma VV              ||      @. \\
0  @>>> \OP3 (-1) @>{\beta '}  >> 4 \OP3 @>{\alpha '}>> \OP3 (1) @>>> 0,
\endCD $$
where  $\gamma = \bmatrix 1&0&0&0\\0&2&0&0\\0&0&1&0\\0&0&0&2
\endbmatrix, \beta ' = \bmatrix -X_1 \\ X_0\\-X_3\\X_2 \endbmatrix,  \alpha ' =
\bmatrix X_0 & X_1 & X_2 & X_3 \endbmatrix.$ This new monad has good reduction
modulo 2, and thus \E\ has a good monad reduction modulo 2.
\medn
3.4)   Consider the bundle \E\ defined over $\Bbb Q$ by the monad
 $$ 0 @>>> \OP3 (-1) @> (-2X_1,2X_0,-X_3,X_2)^\vee >>  4 \OP3 @> (X_0,X_1,X_2,X_3
)>> \OP3 (1) @>>> 0. $$ 
I claim that this \E\ does not have a good monad 
reduction
modulo 2.
\parn
\underbar{Proof}: Of course, this particular monad is also a lift to $\Bbb Z$
which does not reduce well modulo 2. However,
it may be that some lift of an equivalent monad  may exist, that reduces well. So suppose that $A$
is some integral domain between $\Bbb Z$ and $\Bbb Q$ and suppose there is a
lift to
\pew 3A of a monad for \E\ which reduces well modulo 2. This lift
will have matrices $\beta_A$ and $\alpha_A$. The two monads are equivalent over
$\Bbb Q$, hence $\alpha_A = [X_0,X_1,X_2,X_3]\psi$ for some matrix $\psi$ in
$GL(4,\Bbb Q)$. Clearly $\psi$ can be found with entries in $A$. Since the
right hand maps of both monads are surjective modulo 2, the determinant of
$\psi$ is non-zero modulo 2.
By localizing at the multiplicative set obtained from the determinant of $\psi$,
we may assume that $\psi$ is invertible over $A$ and 2 is still in
Spec($A$). Then we see that $c\beta_A = \psi^{-1} [-2X_1,2X_0,-X_3,X_2]^\vee$ 
where
$c \in \Bbb Q$. Write $c$ as $\frac ef$, a ratio of integers in lowest terms.
If 2 divides $e$, then $[-2X_1,2X_0,-X_3,X_2]^\vee$ is identically zero modulo 2,
which is not true. So $e$ is invertible at 2. Hence modulo 2, $\beta_A =
\frac fe \psi^{-1} [0,0,-X_3,X_2]^\vee$. 
This contradicts our assumption that $\beta_A$ modulo 2 is an injection
of bundles. \qed 

\medn
3.5) In the example (3.3), we could also have proceeded as follows: Consider the
automorphism of \pew 3{\Bbb Q} given by $X_1 @>>> X_0, X_1 @>>> X_1/2,
X_2 @>>> X_2, X_3 @>>> X_3/2$. If $\varphi$ is this automorphism, $\varphi^*
\E$ has the monad with $\beta ', \alpha '$ as described there (so that in this case, $\varphi$ fixes \E).   
\medn
3.6) Let \E\ be any  mathematical instanton bundle with $c_1 =0, c_2 = 1$ defined
over a field $k$. Then \E\ has a minimal monad 
$$ 0 @>>> \OP3 (-1) @>\beta >> 4 \OP3 @>\alpha >> \OP3 (1) @>>> 0, $$
where after a change of basis, we may assume that $\alpha = 
\bmatrix X_0 & X_1 & X_2 & X_3 \endbmatrix$. Then the kernel of $\alpha $ is
$\Omega_{\pew 3k}^1(1)$ and the map $\beta $ can be understood as picking out
a section of $H^0(\pew 3k, \Omega_{\pew 3k}^1(2))$. If $V$  is 
the vector space $H^0(\pew 3k, \OP3(1)$ (defined over $k$), $\beta$ can be viewed as picking an element in $\wedge^2V$. The fact that $\beta$ is an inclusion of bundles means that
$\beta$ picks out an indecomposable vector in $\wedge^2V$ (up to scalar multiples). 
Now, the action of $GL(4,k)$ on $\wedge^2V$ has one orbit consisting of the 
indecomposable vectors. Hence, given any \E, there is an automorphism
$\varphi$ of \pew3k such that $\varphi^*\E$ is the standard bundle with the
monad described in (3.3) (ie. with matrix $\beta' = \bmatrix -X_1 \\ X_0\\-X_3\\X_2 \endbmatrix$
.) 
\par
Thus  any \E\ with $c_1 =0, c_2 = 1$ defined over a field of characteristic zero
satisfies the second condition of the Theorem 3.2.
\medn
3.7) A mathematical instanton which satisfies neither of the assumptions of Theorem 3.2.
\par
	Let $k = \Bbb Q(\sqrt{-3})$ and let $A$ be the ring of elements integral over
$\Bbb Z$. The prime number $2$ of $\Bbb Z$ is undecomposed in $A$, hence at the prime $2$, the extension of residue fields has degree 2. Let $l,m$ be the skew lines in \pew 3k, with ideals $(X_0,X_1),
(X_2,X_3)$ respectively. Let $P_1,P_2,P_3,P_4$ be the four points on $l$ with
coordinates $(0,0,1,2)$,$(0,0,1,1)$,$(0,0,1,0)$,\break$(0,0,0,1)$ and $Q_1,Q_2,Q_3,Q_4$ 
the four points on $m$ with coordinates $(1,a,0,0)$,$(1,1,0,0)$,\break$(1,0,0,0)$,
$(0,1,0,0)$, 
where $a \in A$ is an element which reduces modulo 2 to an element $\bar a$ which
is not in $\Bbb Z_2$. The cross-ratio on $l$ of the four points in this order is 
2 and that of the four points  on $m$ is $a$. Recall that the cross-ratio of four distinct 
points on a line is an element of the field which is not 0 or 1 and is invariant
under automorphisms of the line. 
\par
Now let $\phi$ be an automorphism of \pew 3{}
defined over an extension field $k'$ of $k$.  
Let $l',P_i',m',Q_i'$ be the inverse images under $\phi$ of $l,P_i,m,Q_i$. 
The cross-ratio of the appropriate inverse image points is unchanged from
the cross-ratio of the original points. Furthermore, suppose $A'$ is
a ring in $k'$ whose field of fractions is $k'$ and let $\frak p$ be a prime ideal
with residue field of characteristic two. In the light of the remarks
following (3.1), we will assume that $(A', \frak p)$ is a discrete valuation
ring. Then $A'$ contains $A$ and we can find
equations and coordinates for $l',P_i',m',Q_i'$ defined over $A'$. The Hilbert
schemes of lines and points in  \pew 3{A'} are proper over $A'$, hence we can find equations
and coordinates which reduce well modulo $\frak p$ (more concretely, we can 
divide the equations and coordinates by a power of the uniformizing 
parameter of $\frak p$ after which good reduction is possible.) Let $l_0'$,
$m_0'$ be the reductions of $l',m'$ modulo $\frak p$.
On $l_0'$, since the reduced cross-ratio is now zero, the four points are no longer
distinct. On the other hand, the four points on $m'$ do reduce to distinct points
of $m_0'$ by our choice of $a$. 
\par 
 The union $Y = P_1Q_1 \cup P_2Q_2 \cup P_3Q_3 \cup P_4Q_4$ is the zero-scheme of a section $s \in H^0(\pew 3{k}, \E(1))$ for a
mathematical instanton bundle \E\ with $c_2 = 3$, defined over $ k$. In fact,
since $Y$ does not lie on a quadric, \E\ has a unique section up to scalar
multiples.
We will use three facts in the following:
\parn i) The line  $m$ is a jumping line for \E\ of order 3 since it is a
quadri-secant for $Y$. (Indeed, consider the restriction of the sequence for $\I Y$ to $m$:
$ \O m(-1) @>s>> \E_m @>>> \I Y(1)\otimes \O m @>>> 0$.) Hence the same is true
for $m'$ as a jumping line of $\E' = \phi^*\E$. Furthermore, if $\E_{A'}'$ is
a vector bundle on \pew 3{A'}, then the specialization $m_0'$ of $m'$ after
reducing modulo $\frak p$ will be a jumping line for $\E_0'$ of order at least
3. 
\parn ii) Let $s$ be a section in degree one for a mathematical instanton \F\ 
and let $Y$ be the  zero-scheme. Suppose $C$ is a reduced subscheme of $Y$ which
lies on a plane. Then $C$ is a line. For if $H$ is the plane so obtained, when
$s$ is restricted to a section of $\F_H$, it vanishes along the curve $C$, hence
is divisible by the equation of $C$. So $\F_H(-d)$ has a section where $d$ is 
the degree of $C$. The result now follows from the restriction sequence of \F\
to $H$. 
\parn iii) 
A result of  N\"u\ss ler and Trautmann  (true in any characteristic) states that if \F\
is a mathematical instanton with a section $s$ in degree one and if $m$ is any
line contained in the support of the zero scheme of $s$, then $\F_m = \O m(-1)
\oplus \O m(1)$ [N-T].
\par
\underbar{Claim:}
	Let $k'$ be an extension field of $k$, and let $\phi$ be any 
automorphism of  \pew 3{k'}. Then $\E'= \phi^*\E$ does not have a good monad
reduction modulo two. 
\par
\underbar{Proof:}
	Assume the contrary. 
By this assumption, using the same notation as above,  there is
a vector bundle $\E_{A'}'$ on  \pew 3{A'} which specializes to a mathematical
instanton $\E_0'$ on $\pew 3{k(\frak p)}$. Along with a lift of the monad to
$A'$, we can lift the section $s'$ of $\E'(1)$ to $A'$. Hence the specialization of $s'$
will give a section $s_0'$ of  $\E_0'(1)$ which defines a codimension two
zero-scheme $Y_0'$
in $\pew 3{k(\frak p)}$. By our choice of cross-ratios, the lines in $Y_{A'}$
(the zero-scheme of $s'$) 
cannot specialize to four disjoint lines in $\pew 3{k(\frak p)}$. In fact, we claim that
$m_0'$, the specialization of $m'$ is contained in $Y_0'$.
\par
	Indeed, the four points $Q_1',Q_2',Q_3',Q_4'$ on $m'$ reduce to four distinct 
points on $m_0'$, while the four points on $l'$ don't. Since $Y_0'$ cannot have
a component whose reduced subscheme is planar of degree $\geq 2$, in particular
$Y_0'$ cannot contain two distinct lines meeting at a point. So the only 
conclusion is that $m_0'$ lies inside $Y_0'$. But this is a contradiction, since
on the one hand the restriction of $\E_0'$ to $m_0'$ splits as $\O {}(-1) \oplus
\O {}(1)$ using (iii) above and on the other hand
as  $\O {}(-3)\oplus  \O {}(3)$ since the limit of a jumping line of order 
3
is a jumping line of order at least three (and hence equal to 3 as the maximal
possible order of a jumping line is 3.) 
\par
Thus the claim. \qed

\medskip
Lastly, we get the following consequence of Corollary 2.8.
\medn
3.8)\underbar{Theorem}: Let $k$ be a field of characteristic zero. Suppose N  is a component of $\M_{\pew 3k}(c_1,n)$ which contains one bundle \E\ which
reduces modulo two to a stable bundle. Then $N$ has dimension bounded above
by $n^2+4n+1$ for $c_1=0$ and by $n^2 +4n -1$ for $c_1 = -1$.
\parn
\underbar{Proof}: Our definition of good reduction here is more general than
the definition in (3.1). It means merely that there is a vector bundle $\E_A$
on \pew 3A such that $\E_A \otimes _A k$ equals \E, and such that $\E_A \otimes
_{A} k(\frak p)$ is stable. The bounds for $h^1(\Hom (\E,\E))$ of (2.8) are
valid over $k$ as well, by upper semi-continuity. \qed

\medn
3.9)\underbar{Remarks}: We do not know if there are components $N$ which violate
the condition of Theorem 3.8. The condition of degenerating to a stable bundle
can be relaxed to one of degenerating to a semi-stable bundle without any great
change in the dimension bound, since Theorem 2.7 still applies. A 
remark similar to the above theorem can be made about mathematical instantons. Let $N$ be a component of
the moduli space of mathematical instantons with $c_1 =0, c_2=n$ in characteristic
zero. Suppose
that $N$ contains one bundle \E\ which has a good monad reduction modulo
two. Then $N$ is generically smooth of dimension $8n-3$. 
Of course this is well known for the usual component of instantons, ie. the
component containing the bundles corresponding to skew lines. Even in the 
example we gave in (3.7), if the example is deformed in moduli by changing
the cross-ratio of 2 to a value not in $\Bbb Z_2$, we end up with
a bundle with good monad reduction. It may well
be (as seems to be generally expected) that there is only one component for
this moduli space, in which case this remark gives nothing new.
\medn
\underbar{{\bf References:}}
\parn
[A-O] {Ancona,V\. and Ottaviani,G\.,} {\sl On singularities of ${\Cal M}\sb {
{\bold  P}\sp 3}(c_1,c_2)$,\/} preprint, (1995), alg-geom/9502008.

\parn
[B-E]  Barth,W\. and Elencjwag,G\., {\sl Concernant la cohomologie des
        fibr\'es alg\'ebriques stables sur $\bold P_n(\bold C)$,\/} Springer
        Lecture Notes 683 (1978), 1-24.
 
\parn
[B-H]  Barth, W.\ and Hulek, K., {\sl Monads and moduli of vector
bundles}, manuscripta math.\ 25, (1978), 323-347.
 
\parn
[E] {Ein,L.,} {\sl Generalised null correlation
bundles,\/} Nagoya Math\. J\. {\bf 111} (1988), 13--24.

\parn
[H-1]   {Hartshorne,R.,} {\sl Algebraic Geometry,\/}
Springer-Verlag, New York (1977)

\parn
[H-2]   {Hartshorne,R.,} {\sl Stable vector bundles of rank 2
on \pee3,\/} Math\. Ann\.  {\bf 238} (1978), 229--280.

\parn
[H-3]   {Hartshorne,R.,} {\sl Stable reflexive sheaves,\/}
Math\. Ann\. {\bf 254} (1980), 121--176.

\parn
[H-R]  {Hartshorne, R. and Rao,A\.P\.,}  {\sl Spectra and monads of stable bundles,\/}
J\. Math\. Kyoto U\., {\bf 31} (1991), 789--806.

\parn
[O-S] {Okonek,C\. and Spindler,H\.,} {\sl Mathematical instanton bundles on \pee
{2n+1},\/} J\. Reine Angew\. Math. {\bf 364} (1986), 35--50.

\parn
[O-S-S] {Okonek,C\., Schneider,M\. and Spindler,H\.,} {\sl Vector bundles on 
complex projective spaces,\/} Birkh\"auser, Boston and Basel (1980).

\parn
[R-1] {Rao,A\.P\.,} {\sl A note on cohomology modules of rank two bundles,\/} J. of Alg.,
{\bf 86} (1984), 23--34.

\parn
[R-2] {Rao,A\.P\.,} {\sl Mathematical instantons with maximal order jumping
lines,\/} Pac\. J\. Math\., {\bf 178} (1997), 331--344. 
\parn
[L]  {Le Potier,J\.,} {\sl Sur l'espace de modules des fibr\'es
de Yang et Mills,\/} Progress in Math\. {\bf 37}, Birkh\"auser (1983), 65--137.

\parn
[Ma]  {Maggesi,M\.,} {\sl ${\Cal M}\sb {{\bold P}\sp 3}(0,2d\sp 2)$ is 
singular,\/} Forum Math. {\bf 8} (1996), 397--400.

\parn
[M]  {Maruyama, M\.,} {\sl Moduli of stable sheaves I,\/} J\. Math\. Kyoto Univ\.
{\bf 17} (1977), 91--126. 
\parn
[N-T]  {N\"u\ss ler,T\. and  Trautmann,G\.,} {\sl Multiple Koszul
structures on lines and instanton bundles,\/} Int\.J\.Math\., {\bf 5} (1994),
373--388.

\end